\documentclass[11pt,leqno,twoside]{article}
\usepackage{amsthm,amsmath,amssymb,mathrsfs,bbold}
\usepackage[all]{xy}

\title{A Simple Condition for Bounded Displacement}
\author{
Yaar Solomon 
	\vspace{0.1cm}\\
Department of Mathematics,\\
Ben-Gurion University of The Negev,\\
Beer-Sheva, Israel.
	\vspace{0.1cm}\\
yaar.solomon@gmail.com
}

\newenvironment{diagram*}[1]{\begin{align*}{\xymatrix{#1}} \end{align*}} {}

\newenvironment{diagramtow*}[2]{\begin{align*}{\xymatrix{#1}} \qquad {\xymatrix{#2}} \end{align*}} {}

\newcommand{\ceil}[1]{\left\lceil{#1}\right\rceil}

\newcommand{\norm}[1]{\left\|{#1}\right\|}
\newcommand{\inpro}[2]{\langle{#1},{#2}\rangle}

\newcommand{\absolute}[1] {\left|{#1}\right|}
\newcommand{\N}{{\mathbb{N}}}
\newcommand{\Z}{{\mathbb{Z}}}

\newcommand{\R}{{\mathbb{R}}}

\newcommand {\ignore}[1]  {}

\newtheorem{thm}{Theorem}[section]
\theoremstyle{definition}
\newtheorem{definition}[thm]{Definition}

\theoremstyle{plain}
\newtheorem{lem}[thm]{Lemma}

\newtheorem{prop}[thm]{Proposition}
\newtheorem{cor}[thm]{Corollary}

\newtheorem{example}[thm]{Example}
\newtheorem{remark}[thm]{Remark}
\newtheorem*{Remark}{Remark}

\newtheorem{question}[thm]{Question}

\theoremstyle{remark}

\begin{document}
\maketitle

\begin{abstract}
We study separated nets $Y$ that come from primitive substitution tilings of the Euclidean space $\R^d$. We show that the question whether $Y$ is a bounded displacement of $\Z^d$ or not can be reduced, in most cases, to a simple question on the eigenvalues and eigenspaces of the substitution matrix. 
\end{abstract}

\section{Introduction}\label{sec:Introduction}
We denote by $\R^d$ the $d$-dimensional Euclidean space and by $\Z^d$ integer lattice in it. $d(\cdot,\cdot)$ denotes the standard Euclidean metric and $B(x,r)$ is the ball of radius $r$ around $x$, with respect to that metric. We also denote by $\mu_s(\cdot)$ the $s$-dimensional Lebesgue measure in $\R^d$.

A subset $Y\subseteq\R^d$ is a \emph{separated net} if it is uniformly discrete and relatively dense. That is, there exist constants $r,R>0$ such that for any $y_1,y_2\in Y$ we have $d(y_1,y_2)\ge r$, and for every $x\in\R^d$ we have $d(x,Y)\le R$. We say that $Y_1$ is a \emph{bounded displacement (BD)} of $Y_2$ if there is a constant $\alpha$ and a bijection $\phi:Y_1\to\alpha\cdot Y_2$ such that $\sup_{y\in Y_1}\left\{d(y,\phi(y))\right\}<\infty$. 

This paper deals with the following question:  
\begin{question}\label{ques:main_question}
Given a separated net $Y\subseteq\R^d$, is there a BD between $Y$ and $\Z^d$? 
\end{question}
The motivation for this question comes originally from a related question that was asked by Gromov: Is every separated net $Y\subseteq\R^d$ biLipschitz equivalent to $\Z^d$? This question was answered negatively (for $d>1$) by Burago and Kleiner in \cite{BK98}, and independently by McMullen in \cite{McM98}. When considering separated nets, BD equivalence implies biLipschitz equivalence, and this implies that there exists separated net in $\R^d$ which are not BD of $\Z^d$.

In the context of the above questions, it is equivalent to consider tilings of $\R^d$ with finitely many tiles, up to isometry. Obviously, a tiling $\tau$ of $\R^d$ gives rise to separated nets $Y_\tau$ by placing a point in each tile (up to BD). On the other hand, a separated net defines a tiling of $\R^d$ by taking the Voronoi cells. A similar argument gives a tiling with finitely many tiles: Divide the plane to small enough dyadic cubes $Q$. For every $y\in Y$ assign the tile 
\[T_y=\bigcup\left\{\mbox{cubes } Q: Q\mbox{ is closer to } y \mbox{ than to any other } z\in Y\right\}.\]
Denote this tiling by $\tau_Y$, then it is easy to see that any separated net $Y_{\tau_Y}$ is a BD of $Y$. 

We restrict ourselves to a subset of tilings - substitution tilings (see \S\ref{sec:Definitions}). In this context, Theorem \ref{thm:Goal} answers Question \ref{ques:main_question} almost completely. Substitution tilings has a corresponding matrix, the substitution matrix, which we denote by $A_H$ (see Definition \ref{def:SubMatrix+Primitive}). We denote by $\lambda_1,\ldots,\lambda_n$ the eigenvalues of $A_H$, with a descending order in absolute value. These parameters play an important role in our main results, and in the previous related results. 

Question \ref{ques:main_question} was previously studied in the context of substitution tilings in \cite{S11} and \cite{ACG11}. It was shown in \cite{S11} that any primitive substitution tiling, with a matrix $A_H$ of Pisot type, gives rise to a separated net which is a BD of $\Z^d$. Recently Aliste-Prieto, Coronel and Gambaudo have improved this result. They showed that the same holds if $\absolute{\lambda_2}<\lambda_1^{1/d}$, see \cite{ACG11} (note that $\lambda_1>1$). Our main goal is to prove Theorem \ref{thm:Goal}, which extends the results of \cite{S11} and \cite{ACG11} to a wider class, and gives the tight inequality on the eigenvalues for when BD to $\Z^d$ exists. 

We denote by $W_\lambda$ the eigenspace that corresponds to $\lambda$, by $W^\perp$ the subspace which is orthogonal to $W$ with respect to the standard inner product $\inpro{\cdot}{\cdot}$, and let $\mathbb{1}=\begin{pmatrix}1\\ \vdots\\1\end{pmatrix}\in\R^d$.
\begin{thm}\label{thm:Goal}
For a primitive substitution tiling of $\R^d$, fix $t\ge 2$ to be the minimal index that satisfies $W_{\lambda_t}\nsubseteq\mathbb{1}^\perp$. Then the corresponding separated net $Y$ satisfies the following: 
\begin{itemize}
\item[(I)]
If $\absolute{\lambda_t}>\lambda_1^{\frac{d-1}{d}}$ then $Y$ is not a BD of $\Z^d$.
\item[(II)]
If $\absolute{\lambda_t}<\lambda_1^{\frac{d-1}{d}}$ then $Y$ is a BD of $\Z^d$. 
\item[(III)]
If $\absolute{\lambda_t}=\lambda_1^{\frac{d-1}{d}}$ and $\lambda_t$ has a non-trivial Jordan block with at least two (generalized) eigenvectors not in $\mathbb{1}^\perp$, then $Y$ is not a BD of $\Z^d$. Moreover, there are cases where the same consequence holds, and $\lambda_t$ has a trivial Jordan block.

\end{itemize} 
\end{thm}
\begin{remark}\label{remark:Remarks_on_the_main_Thm}
\begin{itemize}
\item
Note that $t=2$ for almost every matrix $A_H$. 
\item
It is follows from the proof of $(II)$ that if there is no $t$ as above, namely $W_{\lambda_t}\subseteq\mathbb{1}^\perp$ for every $t\neq 1$, then $Y$ is a BD of $\Z^d$.
\item
In the case of equality $\absolute{\lambda_t}=\lambda_1^{\frac{d-1}{d}}$, we do not know if there is an example in which $Y$ is a BD of $\Z^d$. 
\end{itemize}
\end{remark}
The proof of the theorem rely on the following result of Laczkovich:
\begin{thm}[\cite{L92}]\label{thm:Laczkovich_Thm}
For a separated net $Y\subseteq\R^d$ and $\beta>0$ the following statements are equivalent:
\begin{itemize}
\item[(i)] There is a constant $C$ such that for any measurable set $A\subseteq\R^d$ we have  
\[\absolute{\#({Y\cap A})-\beta\cdot\mu_d(A)}\le C\cdot\mu_d\left(\{x\in\R^d:d(x,\partial A)\le 1\}\right).\]
\item[(ii)] There is a constant $C$ such that for every finite union of unit lattice cubes $U$ we have  
\[\absolute{\#({Y\cap U})-\beta\cdot\mu_d(U)}\le C\cdot\mu_{d-1}(\partial U).\]
\item[(iii)] There is a BD $\phi:Y\to\beta^{-1/d}\Z^d$.
\end{itemize}
\end{thm}

The organization of this paper is as follows: In section \S\ref{sec:Definitions} we recall the relevant definitions and a few results on substitution tilings. In \S\ref{sec:Economic_Packing_for_Patches} we get a series of estimates that are needed for the proof of Theorem \ref{thm:Goal}. Among them, we prove an isoperimetric lemma, and then use it to generalize a result of Laczkovich. In \S\ref{sec:Proofs_Main_Theorems} we prove Theorem \ref{thm:Goal}, and examples for the different cases this Theorem are given in \S\ref{sec:Examples}.

\medskip{\bf Acknowledgements:} This research was supported by the Israel Science Foundation, grant $190/08$. I wish to thank my supervisor, Barak Weiss, for his many helpful remarks and ideas. I wish to thank Eli Shamovic and Roi Livni for helpful conversations. I also wish to thank David Freeman, who drawn my attention to a mistake in the previous version (the one that was published), that has been corrected here\footnote{The main corrections are in the notations subsection in \S\ref{sec:Definitions}, in Lemma \ref{lem:Discrepancy_Estimate}, in section (III) of the Theorem \ref{thm:Goal}, and in Example \ref{example2}. I also added another example, Example \ref{example4}}.

\section{Preliminaries}\label{sec:Definitions}
A set $T\subseteq\R^d$ is a \emph{tile} if it is biLipschitz homeomorphic to a closed $d$-dimensional ball. Note that this requirement implies in particular that a tile's boundary has a well defined $d-1$-dimensional volume. A \emph{tiling} of a set $U\subseteq\R^d$ is a countable collection of tiles, with pairwise disjoint interiors, such that their union is equal to $U$. 
A tiling $P$ of a bounded set $U\subset\R^d$ is called a \emph{patch}. We call the set $U$ the \emph{support of $P$} and we denote it by $supp(P)$. 
Given a collection of tiles $\mathcal{F}$, we denote by $\mathcal{F}^*$ the set of all patches by the elements of $\mathcal{F}$. For further reading on tiling see for instance \cite{GS87}.

\subsection*{Substitution Tilings}
Let $\xi>1$ and let $\mathcal{F}=\{T_1,\ldots,T_k\}$ be a set of $d$-dimensional tiles. 
\begin{definition}\label{def:SubRule}
A \emph{substitution} is a mapping $H:\mathcal{F}\to\xi^{-1}\mathcal{F}^*$ such that $supp(T_i)=supp(H(T_i))$ for every $i$. Namely, it is a set of dissection rules that shows us how to divide the tiles to other tiles from $\mathcal{F}$ with a smaller scale. We also allow to apply $H$ to finite or infinite collections of tiles. The constant $\xi$ is called the \emph{inflation constant} of $H$.
\end{definition}

\begin{definition}\label{def:SubsTiling}
Let $H$ be a substitution defined on $\mathcal{F}$. Consider the following set of patches: \[\mathcal{P}=\left\{(\xi H)^m(T):m\in\N\:,\: T\in\mathcal{F}\right\}.\] The \emph{substitution tiling space} $X_H$ is the set of all tilings of $\R^d$ that for every patch $P$ in them there is a patch $P'\in\mathcal{P}$ such that $P$ is a sub-patch of $P'$. Every tiling $\tau\in X_H$ is called a \emph{substitution tiling} of $H$. 
\end{definition}

Consider the following equivalence relation on tiles: $T_i\sim T_j$ if there exists an isometry $O$ such that $T_i=O(T_j)$ and $H(T_i)=O(H(T_j))$. We call the representatives of the equivalence classes \emph{basic tiles}, and denote them by $\{\mathcal{T}_1,\ldots,\mathcal{T}_n\}$. By this definition, we can also think of $H$ as a dissection rule on the basic tiles and extend it to collections of tiles as before. For a tile $T$ in the tiling we say that $T$ is \emph{of type $i$} if it is equivalent to $\mathcal{T}_i$. 

\subsection*{Matrices of Substitution} 
\begin{definition}\label{def:SubMatrix+Primitive}
Let $\mathcal{F}=\{\mathcal{T}_1,\ldots,\mathcal{T}_n\}$ be the set of basic tiles. Define the \emph{substitution matrix} of $H$ to be an $n\times n$ matrix, $A_H=(a_{ij})$, where $a_{ij}$ is the number of basic tiles in $\xi H(\mathcal{T}_j)$ which are of type $i$. We say that $H$ is \emph{primitive} if $A_H$ is primitive. That is, if there exists an $m\in \N$ such that $A_H^m>0$.
\end{definition}   

Denote by $e_i$ the $i$'th element of the standard basis of $\R^n$. We use vectors to represent the number of basic tiles from each type in a given patch. For instance, $e_i$ represents one tile of type $i$. Then $A_H(e_i)$ is the $i$'th column of $A_H$. Thus, multiplying the vector $e_i$ by $A_H$ gives a vector that represents the number of basic tiles of each type obtained after applying $H$ on $\mathcal{T}_i$. By linearity, this idea extends to any vector in $\R^n$.

\subsection*{Notations and Previous Results} 
A substitution tiling has many parameters that we need throughout the proofs. For the convenience of the reader we assemble all the notations regarding the parameters of the tiling here.

Our given tiling is denoted by $\tau$ or $\tau_0$, and we fix a separated net $Y$ that correspond to $\tau$. The basic tiles are $\mathcal{F}=\{\mathcal{T}_1,\ldots,\mathcal{T}_n\}$, and $s_1,\ldots,s_n$ denotes their $d$-dimensional volume. $H$ is the substitution, which is always assumed to be primitive, and $\xi>1$ is the inflation constant. We denote by $\lambda_1,\ldots,\lambda_n$ the eigenvalues of $A_H$ in a descending order in absolute value. It is easy to see that $\lambda_1=\xi^d>1$. It also follows from the Perron Frobenius Theorem that $\lambda_1$ is of multiplicity one, and it has positive eigenvector $v_1$ (see \S $3$ in \cite{S11} for details). We fix a Jordan basis of $A_H$ and denote by $v_i$ the $i$'th vector in it ($v_i$ corresponds to $\lambda_i$), and by $v_i(j)$ its $j$'th coordinate. Without loss of generality $v_1(1)=1$. Denote by $u_1=\begin{pmatrix}s_1\\ \vdots\\ s_n\end{pmatrix}$, then it is easy to see that $u_1$ is the left eigenvector of $A_H$ that corresponds to $\lambda_1$. For each $i\in\{1,\ldots,n\}$ we denote by 
\begin{equation}\label{eq:k_i}
k_i= \text{the length of the Jordan chain of } v_i, \text{ counting from } v_i. 
\end{equation}
\begin{Remark}
In our notations the {\bf first} vector in a Jordan chain is the eigenvector $w_1$, the second is a vector $w_2$ that satisfies $(A-\lambda I)w_2=w_1$, and so on. 
\end{Remark} 
Finally, we fix 
\begin{equation}\label{eq:alpha}
\alpha=\frac{\sum_{i=1}^nv_1(i)}{\sum_{i=1}^nv_1(i)\cdot s_i}=\frac{\inpro{\mathbb{1}}{v_1}}{\inpro{u_1}{v_1}}.
\end{equation}
This $\alpha$ is the asymptotic density of $Y$.

\begin{prop}\label{prop:tau_m_Def}
If $H$ is a primitive substitution then $X_H\neq\emptyset$ and for every $\tau\in X_H$ and for every $m\in\N$ there exists a tiling $\tau_m\in X_H$ that satisfies $(\xi H)^m(\tau_m)=\tau$.   
\end{prop} 
\begin{proof}
See \cite{Ro04}.
\end{proof}

Given a tiling $\tau=\tau_0\in X_H$, for every $m$ we fix a tiling $\tau_m$ as in Proposition \ref{prop:tau_m_Def}. $\mathscr{T}^{(m)}$ denotes the set of all tiles of $\tau_m$, and $\mathscr{T}=\bigcup_m\mathscr{T}^{(m)}$. The set of all finite unions of tiles of $\tau_0$ is denoted by $\mathscr{V}$.

We prove Theorem \ref{thm:Goal} using Theorem \ref{thm:Laczkovich_Thm}. To use it we need to estimate the discrepancy $\absolute{\#({Y\cap U})-\alpha\cdot\mu_d(U)}$ for different sets $U$. Note that for every patch $V\in\mathscr{V}$ we have \begin{equation}\label{eq:a_V}
\#(Y\cap V)=\sum_{i=1}^na_i=\inpro{\mathbb{1}}{a_V},\quad\mbox{and}\quad \mu_d(V)=\sum_{i=1}^na_i\cdot s_i=\inpro{u_1}{a_V},
\end{equation}
where $a_V=\begin{pmatrix}a_1\\ \vdots\\a_n\end{pmatrix}$, and $a_j$ is the number of tiles of $\tau_0$ from type $j$ in $V$. Then the discrepancy of $V$ depend only on $a_V$, and is given by the absolute value of the following linear functional:
\begin{equation}
disc(a_V)=\inpro{\mathbb{1}}{a_V}- \frac{\inpro{\mathbb{1}}{v_1}}{\inpro{u_1}{v_1}}\inpro{u_1}{a_V}.
\end{equation}

\begin{lem}\label{lem:Discrepancy_Estimate}
Given a primitive substitution $H$ on $n$ tiles, with a substitution matrix $A_H$, let $t\ge 2$ be the minimal index such that $W_{\lambda_t}\nsubseteq\mathbb{1}^\perp$, and let 
\begin{equation}\label{eq:k}
k=\max\{k_i:v_i\in W_{\lambda_t}, \text{ and } v_i\notin\mathbb{1}^\perp\},
\end{equation} 
where $k_i$ are as in \ref{eq:k_i}. Then there are constants $A_1,A_2>0$, depending only on the parameters of the tiling, with the following properties:
\begin{itemize}
\item[(i)] 
There exists a $j\in\{1,\ldots,n\}$ such that for every $m$ and $T\in\mathscr{T}^{(m)}$ of type $j$
\begin{equation}\label{eq:Improved_Lemma_lower_bound}
A_1\cdot m^{k-1}\absolute{\lambda_t}^m\le\absolute{\#(Y\cap T)-\alpha\cdot\mu_d(T)},
\end{equation}
\item[(ii)]
For every $T\in\mathscr{T}^{(m)}$
\begin{equation}\label{eq:Improved_Lemma_upper_bound}
\absolute{\#(Y\cap T)-\alpha\cdot\mu_d(T)}\le A_2\cdot m^{k-1}\absolute{\lambda_t}^m,
\end{equation} 
\end{itemize}
\end{lem} 
\begin{proof}
Let $T\in\mathscr{T}^{(m)}$ and write $a_T=\sum_{i=1}^nc_iv_i$. 
Note that $disc(v_1)=0$, and also $\inpro{u_1}{v_i}=0$ for every $i\neq 1$. So we have
\begin{equation}\label{eq:disc(a_T)}
disc(a_T)=\inpro{\mathbb{1}}{\sum_{i=2}^nc_iv_i}= \inpro{\mathbb{1}}{\sum_{i=t}^nc_iv_i}.
\end{equation} 
But, if $T$ in $\tau_m$ is of type $j$ then $a_T=A_H^me_j$. Write $e_j=\sum_{i=1}^nb_iv_i$, then \[a_T=A_H^m\left(\sum_{i=1}^nb_iv_i\right)=\sum_{i=1}^nb_iA_H^m(v_i).\]
So for every $i$, $c_i=Const\cdot m^{k_i-1}\cdot\lambda_i^m$ (where the constant on the right hand side also contains a combination of the $b_i$'s). Considering (\ref{eq:disc(a_T)}), this proves $(ii)$. For $(i)$, note that $k=k_\ell$ for some index $\ell$. Since $\{v_1,\ldots,v_n\}$ is a basis of $\R^n$, there exists a $j$ with $b_\ell^{(j)}\neq 0$ in the presentations $e_j=\sum_{i=1}^nb_i^{(j)}v_i$. Using (\ref{eq:disc(a_T)}) in the same way again, we deduce $(i)$ for that particular index $j$.
\end{proof}
\begin{remark}\label{remark:t_does_not_exist}
By (\ref{eq:disc(a_T)}), if $t$ as above does not exist, then the lemma holds with $\lambda_t=0$.
\end{remark}

\section{Economic Packing for Patches}\label{sec:Economic_Packing_for_Patches}
We denote by $\partial A$ and $int(A)$ the boundary and interior of a set $A\subseteq\R^d$ (in the standard topology of $\R^d$), and by $\norm{\cdot}_1$ the standard $\ell_1$ norm on $\R^d$. 

In this section we prove a series of lemmas that help us estimate the terms that appears in Theorem \ref{thm:Laczkovich_Thm}. Our main objective of this section is to prove Proposition \ref{prop:Lemma_3.2_tiles} below. This Proposition gives a very good bound for the number of tiles from each level that one needs in order to obtain a given patch in a substitution tiling. Laczkovich proved this Proposition for the lattice unit cube tiling in \cite{L92}, and here we give a proof for the more general case by generalizing his arguments to our context. Proposition \ref{prop:Lemma_3.2_tiles} is the key point for the proof of Theorem \ref{thm:Goal} in \S\ref{sec:Proofs_Main_Theorems}. 
\begin{lem}\label{lem:Laczkovich_2.1+2.2}
For every $d$ there is a constant $C_1$ such that for every $U\subseteq\R^d$, a finite union of lattice unit cubes, and every $R>0$, we have
\[\mu_d\left(\{x\in U:d(x,\partial U)\le R\}\right)\le C_1\cdot R^d\cdot\mu_{d-1}(\partial U).\]
\end{lem}
\begin{proof}
This is a direct consequence of Lemmas $2.1$ and $2.2$ of \cite{L92}.
\end{proof}
\begin{lem}\label{lem:Distance_1_from_boundary_estimate}
There is a constant $C_2$, that depends on the parameters of the tiling, such that for any $T\in\mathscr{T}$ 
\[\mu_d\left(\{x\in\R^d:d(x,\partial T)\le 1\}\right)\le C_2\cdot\mu_{d-1}(\partial T).\]
\end{lem}
\begin{proof}
Denote by $Q_r$ the $d$-dimensional cube with edge of length $r$. Fix a biLipschitz homeomorphism $\psi_i:\mathcal{T}_i\to Q_1$, denote its biLipschitz constant by $K_i$, and let $K=\max_i\{K_i\}$. Let $T\in\mathscr{T}$ and suppose that $T\in\mathscr{T}^{(m)}$, a tile of type $i$. Then by rescaling the picture by $\xi^m$ we get a biLipschitz homeomorphism $\phi:T\to Q_{\xi^m}$, with the same biLipschitz constant. Since $\phi$ is biLipschitz, it follows that 
\[\phi\left(\{x\in T:d(x,\partial T)\le 1\}\right)\subseteq\{x\in Q_{\xi^m}:d(x,\partial Q_{\xi^m})\le K\}.\]
Then 
\[\mu_d\left(\{x\in T:d(x,\partial T)\le 1\}\right)\le K^d\cdot\mu_d\left(\{x\in Q_{\xi^m}:d(x,\partial Q_{\xi^m})\le K\}\right)\]
Applying the same argument to the tiles which are adjacent to $T$ we obtain 
\begin{equation}\label{eq:1_tube_of_the_boundary_estimate}
\mu_d\left(\{x\in \R^d:d(x,\partial T)\le 1\}\right) \le K^d\cdot\mu_d\left(\{x\in \R^d:d(x,\partial Q_{\xi^m})\le K\}\right).
\end{equation}
It also follows that
\begin{equation}\label{eq:Boundary_of_tile_vs._Boundary_of_cube}
\mu_{d-1}(\partial Q_{\xi^m})\le K^{d-1}\cdot\mu_{d-1}(\partial T)
\end{equation}
(see Theorem $1$ in \cite{EG92} p.$75$). Then by (\ref{eq:1_tube_of_the_boundary_estimate}), (\ref{eq:Boundary_of_tile_vs._Boundary_of_cube}), and Lemma \ref{lem:Laczkovich_2.1+2.2} we have
\[\mu_d\left(\{x\in\R^d:d(x,\partial T)\le 1\}\right) \stackrel{(\ref{eq:1_tube_of_the_boundary_estimate})}\le K^d\cdot
\mu_d\left(\{x\in\R^d:d(x,\partial Q_{\xi^m})\le K\}\right)\le \]
\[C_1\cdot K^{2d}\cdot\mu_{d-1}(\partial Q_{\xi^m}) \stackrel{(\ref{eq:Boundary_of_tile_vs._Boundary_of_cube})}\le 
C_1\cdot K^{3d-1}\cdot\mu_{d-1}(\partial T).\]
\end{proof}


We use the same notations $\mathscr{T}^{(m)}, \mathscr{T}$, and $\mathscr{V}$ as defined at the end of \S\ref{sec:Definitions}.
\begin{lem}\label{lem:Isoperimetric}
Let $T\in\mathscr{T}$, and $c\in(0,\frac{1}{2})$. Then there is an $\varepsilon>0$ such that for any $V\in\mathscr{V}, V\subseteq T$,  with $c\cdot\mu_d(T)\le\mu_d(V)\le\frac{1}{2}\mu_d(T)$, we have
\begin{equation}\label{eq:relative_isoperimetric_for_prototiles}
\mu_{d-1}(\partial V\cap int(T))\ge\varepsilon\cdot\mu_{d-1}(\partial T).
\end{equation}
\end{lem}
\begin{proof}
This Lemma follows from the relative isoperimetric inequality, see \cite{EG92} p.$190$. By this inequality, if $B$ is a closed ball, and $E\subseteq B$ is a closed set of finite perimeter (i.e. $\chi_E$ has a bounded variation) then we have
\begin{equation}\label{eq:relative_isoperimetric_for_ball}
\min\left\{\mu_d(E),\mu_d(B\smallsetminus E)\right\}^{\frac{d-1}{d}}\le C\cdot\mu_{d-1}(\partial E\cap int(B)),
\end{equation}
where $C$ depends only on $d$. Fix a biLipschitz homeomorphism $\psi_i:\mathcal{T}_i\to B(0,1)$, denote its biLipschitz constant by $K_i$, and let $K=\max_i\{K_i\}$. Suppose that $T$ is a tile of type $i$, then by rescaling the picture by $\xi^m$ we get a biLipschitz homeomorphism $\phi:T\to B=B(0,\xi^m)$, with the same biLipschitz constant. Since $\phi$ is biLipschitz, it follows that 
\[\frac{1}{K^d}\mu_d(V)\le\mu_d(\phi(V))\le K^d\mu_d(V)\quad\mbox{and}\quad
\frac{1}{K^{d-1}}\mu_{d-1}(\phi(\partial V\cap int(T)))\le \mu_{d-1}(\partial V\cap int(T))\]
(see \cite{EG92} p.$75$). Considering (\ref{eq:relative_isoperimetric_for_ball}) with $E=\phi(V)$ we obtain
\[\mu_{d-1}(\partial V\cap int(T))\ge
\frac{\min\left\{\mu_d(\phi(V)),\mu_d(\phi(T\smallsetminus V))\right\}^{\frac{d-1}{d}}}{C\cdot K^{d-1}}\ge
\frac{\min\left\{\mu_d(V),\mu_d(T\smallsetminus V)\right\}^{\frac{d-1}{d}}}{C\cdot (K^{d-1})^2}\]
\[\ge\frac{c^{\frac{d-1}{d}}\cdot\mu_d(T)^{\frac{d-1}{d}}}{C\cdot (K^{d-1})^2}= \frac{c^{\frac{d-1}{d}}\cdot s_i^{\frac{d-1}{d}}\cdot\xi^{m(d-1)}}{C\cdot (K^{d-1})^2}=\frac{c^{\frac{d-1}{d}}\cdot s_i^{\frac{d-1}{d}}}{C\cdot (K^{d-1})^2\cdot\mu_{d-1}(\partial \mathcal{T}_i)}\cdot\mu_{d-1}(\partial T).\]
Setting $s=\min_i\{s_i\}$ and $D_{\max}=\max_i\{\mu_{d-1}(\partial \mathcal{T}_i)\}$ we get 
\begin{equation}\label{eq:Epsilon}
\varepsilon=\frac{c^{\frac{d-1}{d}}\cdot s^{\frac{d-1}{d}}}{C\cdot (K^{d-1})^2\cdot D_{\max}}
\end{equation}
that satisfies the assertion, and does not depend on the type of the tile $T$. 
\end{proof}

\begin{cor}\label{cor:Laczkovich_Lemma_3.1}
Let $T\in\mathscr{T}$, $c\in(0,\frac{1}{2})$, and $\varepsilon$ as in (\ref{eq:Epsilon}). Suppose that $V\in\mathscr{V}, V\subseteq T$ with $\mu_d(V)\le (1-c)\cdot\mu_d(T)$ and $\mu_{d-1}(\partial V\cap int(T))<\varepsilon\cdot\mu_{d-1}(\partial T)$, then $\mu_d(V)< \frac{1}{2}\mu_d(T)$.
\end{cor}
\begin{proof}
Assume otherwise, then we have $\mu_{d-1}(\partial (T\smallsetminus V)\cap int(T))<\varepsilon\cdot\mu_{d-1}(\partial T)$ and $c\cdot\mu_d(T)\le\mu_d(T\smallsetminus V)\le\frac{1}{2}\mu_d(T)$, contradicting Lemma \ref{lem:Isoperimetric}.
\end{proof}

For a $T$ in $\tau_m$ we denote by $T^*$ the unique tile of $\tau_{m+1}$ that contains $T$. We denote $\rho=\frac{\max_i\{s_i\}}{\min_i\{s_i\}}\ge 1$, then for any tile $T\in\mathscr{T}$ we have
\begin{equation}\label{eq:T_vs._T^*}
\rho^{-1}\cdot\xi^{-d}\le\frac{\mu_d(T)}{\mu_d(T^*)}\le\rho\cdot\xi^{-d}
\end{equation}
For a set $X\subseteq\mathscr{T}$ we denote by $S(X)$ the closure of $X$ under the operations of disjoint union and proper difference, where every element of $X$ can be used only once. 
For the following lemma we set $\varepsilon$ as in Lemma \ref{lem:Isoperimetric} and define the following constants:
\begin{equation}\label{eq:The_constants}
D_{\min}=\min_i\left\{\mu_{d-1}(\partial \mathcal{T}_i)\right\},\quad
C=\frac{\rho\cdot\xi(\rho\cdot\xi^d+1)}{\varepsilon\cdot D_{\min}}\quad\mbox{and}\quad c=(2\rho)^{-1}\cdot\xi^{-d}\in\left(0,\frac{1}{2}\right).
\end{equation}
\begin{prop}\label{prop:Lemma_3.2_tiles}
Let
\begin{equation}\label{eq:Lemma_assumptions}
V\in\mathscr{V},\quad T\in\mathscr{T},\quad V\subseteq T,\quad\mbox{and}\quad \mu_d(V)\le\frac{1}{2}\mu_d(T).  
\end{equation} 
Then there exists $T_1,\ldots,T_n\in\mathscr{T}$, $T_i\subseteq T$ for all $i$, such that $V\in S(\{T_1,\ldots,T_n\})$, and for every $m$ we have:
\[\#\{i:T_i\in\mathscr{T}^{(m)}\}\le C\cdot
\frac{\mu_{d-1}(\partial V\cap int(T))}{\xi^{m(d-1)}}.\]
\end{prop}
\begin{proof}
The proof is by induction on $m$, where $T\in\mathscr{T}^{(m)}$. If $m=0$ then $\mu_d(V)\le\frac{1}{2}\mu_d(T)$ implies that $V=\varnothing$, so the assertion is obvious. Assume the assertion for any $T\in\mathscr{T}^{(m)}$ with $m<m_0$, and let $V$ and $T$ be as in (\ref{eq:Lemma_assumptions}) with $T\in\mathscr{T}^{(m_0)}$. Consider the following collection of tiles:
\[\mathscr{A}=\left\{P\in\mathscr{T}:\begin{matrix}
P\subseteq T,\\ 
\mu_d(P\cap V)\ge c\cdot\mu_d(P),\\ 
\mu_{d-1}(\partial V\cap int(P))<\varepsilon\cdot\mu_{d-1}(\partial P)
\end{matrix}\right\},\]
where $\varepsilon$ is as in Lemma \ref{lem:Isoperimetric} and $c$ is as in (\ref{eq:The_constants}) (it might be that $\mathscr{A}=\varnothing$). Note that every $P\in\mathscr{A}$ satisfies:
\[\mu_d(P\smallsetminus V)\le (1-c)\mu_d(P), \quad\mbox{and}\quad 
\mu_{d-1}(\partial(P\smallsetminus V)\cap int(P))< \varepsilon\cdot\mu_{d-1}(\partial P).\]
Then by Corollary \ref{cor:Laczkovich_Lemma_3.1} we have
\begin{equation}\label{eq:P_minus_V_estimation}
\mu_d(P\smallsetminus V)<\frac{1}{2}\mu_d(P).
\end{equation}

Let $P_1,\ldots,P_\ell$ be the maximal elements of $\mathscr{A}$ (w.r.t. inclusion). Then $\bigcup\mathscr{A}=\bigcup_{j=1}^\ell P_j\subseteq T$, and $P_1,\ldots,P_\ell$ has pairwise disjoint interiors. Denote $V_1=V\cup\bigcup_{j=1}^\ell P_j$. Then
\begin{equation}\label{eq:V_1_is_strictly_contained_in_T}
\mu_d(V_1)\le\mu_d(V)+\sum_{j=1}^\ell\mu_d(P_j\smallsetminus V)\stackrel{(\ref{eq:P_minus_V_estimation})}< \frac{1}{2}\mu_d(T)+\sum_{j=1}^\ell\frac{1}{2}\mu_d(P_j) \le\mu_d(T).
\end{equation}
Note that if $\mathscr{A}=\varnothing$ we only get $\le$ in the middle inequality, but then the last inequality is strict.
Thus $V_1\subsetneqq T$, and in particular $P_j\subsetneqq T$ for every $j$. By (\ref{eq:P_minus_V_estimation}), we may apply the induction hypothesis for $P_j\smallsetminus V\in\mathscr{V}$ and the tile $P_j$, to obtain tiles $T_{j1},\ldots,T_{jn_j}$ such that $T_{jr}\subseteq P_j$, $P_j\smallsetminus V\in S(\{T_{j1},\ldots,T_{jn_j}\})$, and for every $m$ we have:
\begin{equation}\label{eq:By_the_induc_for_P_j}
\#\left\{r:T_{jr}\in\mathscr{T}^{(m)}\right\}\le 
C\cdot\frac{\mu_{d-1}(\partial V\cap int (P_j))}{\xi^{m(d-1)}}.
\end{equation}  

Now let $T_1,\ldots,T_n$ be the maximal tiles that are contained in $V_1$. Then $T_1,\ldots,T_n$ has pairwise disjoint interiors and their union is equal to $V_1$. So we can write 
\[V=V_1\smallsetminus\bigcup_{j=1}^\ell(P_j\smallsetminus V)= \left(\bigcup_{i=1}^nT_i\right)\smallsetminus\bigcup_{j=1}^\ell(P_j\smallsetminus V),\]
where the sets $P_j\smallsetminus V$ are pairwise disjoint. This implies that 
\[V\in S(\{T_1,\ldots,T_n,T_{11},\ldots,T_{1n_1},\ldots,T_{\ell 1},\ldots,T_{\ell n_\ell}\}).\]

Fix $m\in\N$ and denote $E=\{i:T_i\in\mathscr{T}^{(m)}\}$, and $E_j=\{r:T_{jr}\in\mathscr{T}^{(m)}\}$. It remains to show that 
\begin{equation}\label{eq:goal_of_3.2}
\absolute{E}+\sum_{j=1}^\ell\absolute{E_j}\le C\cdot\frac{\mu_{d-1}(\partial V\cap int (T))}{\xi^{m(d-1)}}.
\end{equation}

We first estimate $\absolute{E}$. Fix an $i\in E$. Since $T_i$ is maximal in $V_1$, if follows that $T_i^*\nsubseteq V_1$. In particular, by the definition of $V_1$, since the $P_j$'s are maximal in $\mathscr{A}$, we have $T_i^*\notin\mathscr{A}$. by (\ref{eq:V_1_is_strictly_contained_in_T}), $V_1\subsetneqq T$, then $T_i\subsetneqq T$, and therefore $T_i^*\subseteq T$. Our next goal is to show that 
\begin{equation}\label{eq:T^*_estimate}
\mu_d(T_i^*\cap V)\ge c\cdot\mu_d(T_i^*).
\end{equation}
If $int(T_i)\cap\left(\bigcup_{j=1}^\ell int(P_j)\right)=\varnothing$ then $T_i\subseteq V$, and therefore 
\[\mu_d(T_i^*\cap V)\ge\mu_d(T_i)\stackrel{(\ref{eq:T_vs._T^*})}\ge \rho^{-1}\cdot\xi^{-d}\cdot\mu_d(T_i^*)\stackrel{(\ref{eq:The_constants})}>c\cdot\mu_d(T_i^*).\]
Otherwise, $int(T_i)$ intersect $int(P_j)$ for some $j$. Then either $T_i\subsetneqq P_j$ or $P_j\subseteq T_i$. If $T_i\subsetneqq P_j$ then $T_i^*\subseteq P_j\subseteq V_1$, a contradiction. Then $P_j\subseteq T_i$ whenever $int(T_i)\cap int(P_j)\neq\varnothing$. Denote by $J$ the set of indices $j$ such that $P_j\subseteq T_i$, then we have
\[\mu_d(T_i\smallsetminus V)\le\mu_d\left(\bigcup_{j\in J}(P_j\smallsetminus V)\right)\le\sum_{j\in J}\mu_d(P_j\smallsetminus V)\]
\[\stackrel{(\ref{eq:P_minus_V_estimation})}<\sum_{j\in J}\frac{1}{2}\mu_d(P_j)\le \frac{1}{2}\mu_d(T_i).\]
Hence
\[\mu_d(T_i^*\cap V)\ge\mu_d(T_i\cap V)>\frac{1}{2}\mu_d(T_i) \stackrel{(\ref{eq:T_vs._T^*})}\ge(2\rho)^{-1}\cdot\xi^{-d}\cdot\mu_d(T_i^*) \stackrel{(\ref{eq:The_constants})}=c\cdot \mu_d(T_i^*).\]
Thus (\ref{eq:T^*_estimate}) holds. Since $T_i^*\subseteq T$ and $T_i^*\notin\mathscr{A}$, it follows form (\ref{eq:T^*_estimate}) and from the definition of $\mathscr{A}$ that $T_i^*\notin\mathscr{A}$ because it satisfies
\begin{equation}\label{eq:T_i^*capV_boundary}
\mu_{d-1}(\partial V\cap int(T_i^*))\ge\varepsilon\cdot\mu_{d-1}(\partial T_i^*).
\end{equation}

Let $K=\partial V\cap\bigcup_{i\in E}int(T_i^*)$. Since the $T_i$'s are distinct elements of $\mathscr{T}^{(m)}$, and by (\ref{eq:T_vs._T^*}), each point of $K$ is covered by at most $\rho\cdot\xi^d$\quad $T_i^*$'s. Therefore, by (\ref{eq:T_i^*capV_boundary}), we have
\[\rho\cdot\xi^d\mu_{d-1}(K)\ge\sum_{i\in E}\mu_{d-1}(K\cap T_i^*)= \sum_{i\in E}\mu_{d-1}(\partial V\cap intT_i^*)
\stackrel{(\ref{eq:T_i^*capV_boundary})}\ge\varepsilon\cdot\mu_{d-1}(\partial T_i^*)\cdot\absolute{E},\]
and hence
\begin{equation}\label{eq:E}
\absolute{E}\le 
\frac{\rho\cdot\xi^d}{\varepsilon\cdot\mu_{d-1}(\partial T_i^*)}\mu_{d-1}(K).
\end{equation}
Now define
\[J_1=\{j:P_j\subseteq T_i^* \mbox{ for some } i\in E\},\quad\mbox{and}\quad J_2=\{1,\ldots,\ell\}\smallsetminus J_1. \]
If $j\in J_1$ and $r\in E_j$ then $T_{jr}\subseteq T_i^*$ for some $i$. Since $T_i^*$ contains at most $\rho\cdot\xi^d$ tiles of $\mathscr{T}^{(m)}$ we have
\[\sum_{j\in J_1}\absolute{E_j}\le\rho\cdot\xi^d\absolute{E}.\]
If $j\in J_2$ and $i\in E$ then $int(P_j)\cap int(T_i^*)=\varnothing$ (since $T_i^*\nsubseteq P_j$). Then the set $K_j=\partial V\cap int(P_j)$ is disjoint from $K$. By (\ref{eq:By_the_induc_for_P_j}) we have $\absolute{E_j}\le C\cdot\frac{\mu_{d-1}(K_j)}{\xi^{m(d-1)}}$, and hence
\begin{equation}\label{eq:The_E_j's}
\sum_{j=1}^\ell\absolute{E_j}=\sum_{j\in J_1}\absolute{E_j}+ \sum_{j\in J_2}\absolute{E_j}\le\rho\cdot\xi^d\absolute{E}+C\cdot \frac{\mu_{d-1}\left(\bigcup_{j\in J_2}K_j\right)}{\xi^{m(d-1)}}.
\end{equation}
The sets $K$ and $\bigcup_{j\in J_2}K_j$ are disjoint, and their union is a subset of $\partial V\cap int(T)$, hence
\[\absolute{E}+\sum_{j=1}^\ell\absolute{E_j}\stackrel{(\ref{eq:The_E_j's})}\le 
(\rho\cdot\xi^d+1)\absolute{E}+\sum_{j\in J_2}\absolute{E_j}\stackrel{(\ref{eq:E}), (\ref{eq:The_E_j's})}\le \frac{\rho\cdot\xi^d(\rho\cdot\xi^d+1)}{\varepsilon\cdot\mu_{d-1}(\partial T_i^*)}\mu_{d-1}(K)+C\cdot \frac{\mu_{d-1}\left(\bigcup_{j\in J_2}K_j\right)}{\xi^{m(d-1)}}\]
\[\stackrel{(\ref{eq:The_constants})}\le\frac{\rho\cdot\xi^d(\rho\cdot\xi^d+1)}{\varepsilon \cdot D_{\min}\cdot\xi^{(m+1)(d-1)}}\mu_{d-1}(K)+
C\cdot \frac{\mu_{d-1}\left(\bigcup_{j\in J_2}K_j\right)}{\xi^{m(d-1)}}\]
\[\stackrel{(\ref{eq:The_constants})}\le\frac{C}{\xi^{m(d-1)}}\left(\mu_{d-1}(K)+ \mu_{d-1}\left(\bigcup_{j\in J_2}K_j\right)\right)
\le C\cdot\frac{\mu_{d-1}(\partial V\cap int(T))}{\xi^{m(d-1)}}.\]
Thus (\ref{eq:goal_of_3.2}) holds and the proof is complete.
\end{proof}

\section{Proof of the Main Theorems}\label{sec:Proofs_Main_Theorems}
\begin{proof}[Proof of Theorem \ref{thm:Goal}]
{\bf\underline{Proof of (I):}}
We show that if $\absolute{\lambda_t}>\lambda_1^{\frac{d-1}{d}}$ then $(i)$ of Theorem \ref{thm:Laczkovich_Thm} does not hold for any $\alpha$. Fix a $j\in\{1,\ldots,n\}$ such that (\ref{eq:Improved_Lemma_lower_bound}) holds for any tile of type $j$, and consider the sequence of measurable sets $T^{(m)}$, the $m$'th inflation of $\mathcal{T}_j$. Then by Lemma \ref{lem:Distance_1_from_boundary_estimate} we have
\[
\mu_d\left(\{x\in\R^d:d(x,\partial T)\le 1\}\right)\le C_2\cdot\mu_{d-1}(\partial T)=C_2\cdot\mu_{d-1}(\partial \mathcal{T}_j)\cdot\xi^{m(d-1)}.
\] 
Recall that $\xi^d=\lambda_1$, then we have
\begin{equation}\label{eq:d-1-dim_vol}
\mu_d\left(\{x\in\R^d:d(x,\partial T)\le 1\}\right)\le C_2\left(\lambda_1^{\frac{d-1}{d}}\right)^m\mu_{d-1}(\partial \mathcal{T}_i).
\end{equation} 
As we did in the proof of Lemma (\ref{lem:Discrepancy_Estimate}), for any $\alpha$ different than the one defined in (\ref{eq:alpha}) $disc(v_1)\neq 0$, and so $disc(a_T)=Const\cdot\lambda_1^m$. For large $m$'s, this is obviously greater than any constant times $\mu_d\left(\{x\in\R^d:d(x,\partial T)\le 1\}\right)$. For $\alpha$ as in (\ref{eq:alpha}), by Lemma \ref{lem:Discrepancy_Estimate} we have
\[\absolute{\#(Y\cap T^{(m)})-\alpha\cdot\mu_d(T^{(m)})}\ge A_1\cdot\absolute{\lambda_t}^m,\]
which by assumption is greater than $Const\left(\lambda_1^{\frac{(d-1)}{d}}\right)^m$, for any constant and for a large enough $m$'s. Considering (\ref{eq:d-1-dim_vol}), we proved that $(i)$ of Theorem \ref{thm:Laczkovich_Thm} does not hold.

{\bf\underline{Proof of (II):}} We show that $(ii)$ of Theorem \ref{thm:Laczkovich_Thm} holds, where $\alpha$ is as in (\ref{eq:alpha}). Let $R=\ceil{\max_i\limits\{diam(\mathcal{T}_i)\}}$, where $diam(A)$ denote the diameter of a set $A$.
It is sufficient to show that $(ii)$ holds for any $U$, a finite union of $R$-cubes (cubes with edge length $R$ and corners at $R\cdot\Z^d$). Let $U$ be a finite union of $R$-cubes. For every $y\in Y$ we denote by $T_y$ the tile of $\tau_0$ that corresponds to $y$, and define an $V\in\mathscr{V}$ by $V=\bigcup\left\{T_y:y\in U\right\}$. Then $U\subseteq V\cup(U\smallsetminus V)$. Note that $U\smallsetminus V\subseteq\{x\in U:d(x,\partial U)\le R\}$, so it follows from Lemma \ref{lem:Laczkovich_2.1+2.2} that
\[\mu_d(U\smallsetminus V)\le 
C_1\cdot R^d\cdot\mu_{d-1}(\partial U).\]
Since $\#(U\cap Y)=\#(V\cap Y)$ we have
\begin{equation}\label{eq:discrepancy_V_vs._discrepancy_U}
\absolute{\#(U\cap Y)-\alpha\cdot\mu_d(U)}\le
\absolute{\#(V\cap Y)-\alpha\cdot\mu_d(V)}+ 
\alpha\cdot C_1\cdot R^d\cdot\mu_{d-1}(\partial U).
\end{equation}
So it is enough to estimate $\absolute{\#(V\cap Y)-\alpha\cdot\mu_d(V)}$.

Next we claim that $\partial V\subseteq\{x\in\R^d:d(x,\partial U)\le R\}$.
Indeed, if $x\in\partial V$ then either $x\in U$ or $x\notin U$. If $x\in U$, since $x\in \partial V$, $x\in\partial T_y$ for some $y\notin U$, and therefore $d(x,\partial U)\le d(x,y)\le diam(T_y)\le R$. A similar argument holds if $x\notin U$ since $x$ also belong to $\partial T_y$ for some $y\in U$. Therefore, every tile $T$ of $\tau_0$ with $T\cap\partial V\neq\varnothing$ is contained in $\{x\in\R^d:d(x,\partial U)\le 2R\}$. Denote by $C_3=\max_i\frac{\mu_{d-1}(\partial\mathcal{T}_i)}{\mu_d(\mathcal{T}_i)}$. Then by Lemma \ref{lem:Laczkovich_2.1+2.2} we have
\begin{equation}\label{eq:partial(V)_vs._partial(U)}
\begin{split}&
\mu_{d-1}(\partial V)\le
\sum_{\substack{T\in\mathscr{T}^{(0)} \\ T\cap\partial V\neq\varnothing}}\mu_{d-1}(\partial T)\le 
\sum_{\substack{T\in\mathscr{T}^{(0)} \\ T\cap\partial V\neq\varnothing}}C_3\cdot\mu_d(T)\le 
\\&
C_3\cdot\mu_d\left(\{x\in\R^d:d(x,\partial U)\le 2R\}\right)\le C_3\cdot C_1\cdot (2R)^d\cdot\mu_{d-1}(\partial U).\end{split}
\end{equation}

To finish the proof, we apply Proposition \ref{prop:Lemma_3.2_tiles} to $V$. We pick a large enough $T\in\mathscr{T}$ such that (\ref{eq:Lemma_assumptions}) holds. By Proposition \ref{prop:Lemma_3.2_tiles} we obtain $T_1,\ldots,T_n\in\mathscr{T}$ such that $V\in S(\{T_1,\ldots,T_n\})$, and for every $m$ we have:
\begin{equation}\label{eq:By_the_Induction_Hypothesis}
\#\{i:T_i\in\mathscr{T}^{(m)}\}\le C\cdot
\frac{\mu_{d-1}(\partial V\cap int(T))}{\xi^{m(d-1)}}.
\end{equation}
Note that if $A,B\in\mathscr{V}$ and $int(A)\cap int(B)=\varnothing$ then
\[\#({Y\cap (A\cup B)})-\alpha\cdot\mu_d(A\cup B)= 
\#({Y\cap A})-\alpha\cdot\mu_d(A)+\#({Y\cap B})-\alpha\cdot\mu_d(B),\]
and similarly if $B\subseteq A$ then 
\[\#({Y\cap (A\smallsetminus B)})-\alpha\cdot\mu_d(A\smallsetminus B)= 
\#({Y\cap A})-\alpha\cdot\mu_d(A)-(\#({Y\cap B})-\alpha\cdot\mu_d(B)).\]
Therefore, since $V\in S(\{T_1,\ldots,T_n\})$, we have 
\[\absolute{\#({Y\cap V})-\alpha\cdot\mu_d(V)}\le\sum_{i=1}^n \absolute{\#({Y\cap T_i})-\alpha\cdot\mu_d(T_i)}\le 
\sum_{m=0}^\infty\sum_{T_i\in\mathscr{T}^{(m)}}\absolute{\#({Y\cap T_i})-\alpha\cdot\mu_d(T_i)}\]
\[\stackrel{(\ref{eq:Improved_Lemma_upper_bound}),(\ref{eq:By_the_Induction_Hypothesis})}\le\sum_{m=0}^\infty \left[C\cdot
\frac{\mu_{d-1}(\partial V\cap int(T))}{\xi^{m(d-1)}}\cdot A_2\cdot m^{k_t-1}\absolute{\lambda_t}^m\right]\le \left[\sum_{m=0}^\infty\frac{m^{k_t-1}\absolute{\lambda_t}^m}{\left(\xi^{d-1} \right)^m}\right]\cdot C\cdot A_2\cdot\mu_{d-1}(\partial V).\]
By the assumption, $\absolute{\lambda_2}<\lambda_1^{\frac{d-1}{d}}=\xi^{d-1}$, and therefore the series converges and we have
\[\absolute{\#({Y\cap V})-\alpha\cdot\mu_d(V)}\le Const\cdot\mu_{d-1}(\partial V).\]
Considering (\ref{eq:discrepancy_V_vs._discrepancy_U}) and (\ref{eq:partial(V)_vs._partial(U)}), we have shown $(ii)$ of Theorem \ref{thm:Laczkovich_Thm}, which implies the assertion.

{\bf\underline{Proof of (III):}} The existence of two (generalized) eigenvectors not in $\mathbb{1}^\perp$ implies that the value of $k$ from (\ref{eq:k}) is at least $2$ (note that having just one eigenvector with $k_i>1$ would suffice). So using (\ref{eq:Improved_Lemma_lower_bound}) in the same way as in the proof of $(I)$ proves the assertion. 

For the case where the Jordan block of $\lambda_t$ is trivial, we give an example in $\R^3$, where the corresponding separated net is not a BD of $\Z^3$.\\
{\bf Example:}
Consider the substitution rule $H$ that is defined by this picture:
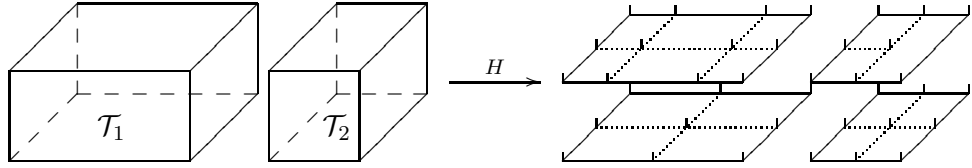
\[\xygraph{
!{<0cm,0cm>;<0.3cm,0cm>:<0cm,0.3cm>::}
!{(-24.5,0)}*{}="1T1"
!{(-16.5,0)}*{}="1T2"
!{(-21.5,3)}*{}="1T3"
!{(-13.5,3)}*{}="1T4"
!{(-24.5,4)}*{}="1T5"
!{(-16.5,4)}*{}="1T6"
!{(-21.5,7)}*{}="1T7"
!{(-13.5,7)}*{}="1T8"
!{(-13,0)}*{}="2T1"
!{(-9,0)}*{}="2T2"
!{(-10,3)}*{}="2T3"
!{(-6,3)}*{}="2T4"
!{(-13,4)}*{}="2T5"
!{(-9,4)}*{}="2T6"
!{(-10,7)}*{}="2T7"
!{(-6,7)}*{}="2T8"
!{(0,0)}*{}="1d1"
!{(4,0)}*{}="1d2"
!{(8,0)}*{}="1d3"
!{(1.5,1.5)}*{}="1d4"
!{(5.5,1.5)}*{}="1d5"
!{(9.5,1.5)}*{}="1d6"
!{(3,3)}*{}="1d7"
!{(7,3)}*{}="1d8"
!{(11,3)}*{}="1d9"
!{(0,3.5)}*{}="1u1"
!{(2,3.5)}*{}="1u2"
!{(6,3.5)}*{}="1u3"
!{(8,3.5)}*{}="1u4"
!{(1.5,5)}*{}="1u5"
!{(3.5,5)}*{}="1u6"
!{(7.5,5)}*{}="1u7"
!{(9.5,5)}*{}="1u8"
!{(3,6.5)}*{}="1u9"
!{(5,6.5)}*{}="1u10"
!{(9,6.5)}*{}="1u11"
!{(11,6.5)}*{}="1u12"
!{(11,0)}*{}="2d1"
!{(13,0)}*{}="2d2"
!{(15,0)}*{}="2d3"
!{(12.5,1.5)}*{}="2d4"
!{(14.5,1.5)}*{}="2d5"
!{(16.5,1.5)}*{}="2d6"
!{(14,3)}*{}="2d7"
!{(16,3)}*{}="2d8"
!{(18,3)}*{}="2d9"
!{(11,3.5)}*{}="2u1"
!{(13,3.5)}*{}="2u2"
!{(15,3.5)}*{}="2u3"
!{(12.5,5)}*{}="2u4"
!{(14.5,5)}*{}="2u5"
!{(14,6.5)}*{}="2u6"
!{(16,6.5)}*{}="2u7"
!{(18,6.5)}*{}="2u8"
!{(0,0.4)}*{}="1d1'"
!{(4,0.4)}*{}="1d2'"
!{(8,0.4)}*{}="1d3'"
!{(1.5,1.9)}*{}="1d4'"
!{(5.5,1.9)}*{}="1d5'"
!{(9.5,1.9)}*{}="1d6'"
!{(3,3.4)}*{}="1d7'"
!{(7,3.4)}*{}="1d8'"
!{(11,3.4)}*{}="1d9'"
!{(0,3.9)}*{}="1u1'"
!{(2,3.9)}*{}="1u2'"
!{(6,3.9)}*{}="1u3'"
!{(8,3.9)}*{}="1u4'"
!{(1.5,5.4)}*{}="1u5'"
!{(3.5,5.4)}*{}="1u6'"
!{(7.5,5.4)}*{}="1u7'"
!{(9.5,5.4)}*{}="1u8'"
!{(3,6.9)}*{}="1u9'"
!{(5,6.9)}*{}="1u10'"
!{(9,6.9)}*{}="1u11'"
!{(11,6.9)}*{}="1u12'"
!{(11,0.4)}*{}="2d1'"
!{(13,0.4)}*{}="2d2'"
!{(15,0.4)}*{}="2d3'"
!{(12.5,1.9)}*{}="2d4'"
!{(14.5,1.9)}*{}="2d5'"
!{(16.5,1.9)}*{}="2d6'"
!{(14,3.4)}*{}="2d7'"
!{(16,3.4)}*{}="2d8'"
!{(18,3.4)}*{}="2d9'"
!{(11,3.9)}*{}="2u1'"
!{(13,3.9)}*{}="2u2'"
!{(15,3.9)}*{}="2u3'"
!{(12.5,5.4)}*{}="2u4'"
!{(14.5,5.4)}*{}="2u5'"
!{(14,6.9)}*{}="2u6'"
!{(16,6.9)}*{}="2u7'"
!{(18,6.9)}*{}="2u8'"
!{(-5,3.5)}*{}="A1"
!{(-1,3.5)}*{}="A2"
!{(-20,1.5)}*{\mathcal{T}_1}="T1"
!{(-10,1.5)}*{\mathcal{T}_2}="T2"
"1d1"-"1d3" "1d4"-@{..}"1d6" "1d7"-"1d9"
"1d1"-"1d7" "1d2"-@{..}"1d8" "1d3"-"1d9"
"1u1"-"1u4" "1u5"-@{..}"1u8" "1u9"-"1u12"
"1u1"-"1u9" "1u2"-@{..}"1u10" "1u3"-@{..}"1u11" "1u4"-"1u12" 
"2d1"-"2d3" "2d4"-@{..}"2d6" "2d7"-"2d9"
"2d1"-"2d7" "2d2"-@{..}"2d8" "2d3"-"2d9"
"2u1"-"2u3" "2u4"-@{..}"2u5" "2u6"-"2u8"
"2u1"-"2u6" "2u2"-@{..}"2u7" "2u3"-"2u8" 
"1T1"-"1T2" "1T1"-@{--}"1T3" "1T2"-"1T4" "1T3"-@{--}"1T4"
"1T5"-"1T6" "1T5"-"1T7" "1T6"-"1T8" "1T7"-"1T8"
"1T1"-"1T5" "1T2"-"1T6" "1T3"-@{--}"1T7" "1T4"-"1T8"
"2T1"-"2T2" "2T1"-@{--}"2T3" "2T2"-"2T4" "2T3"-@{--}"2T4"
"2T5"-"2T6" "2T5"-"2T7" "2T6"-"2T8" "2T7"-"2T8"
"2T1"-"2T5" "2T2"-"2T6" "2T3"-@{--}"2T7" "2T4"-"2T8"
"A1":^H"A2"
"1d1"-"1d1'" "1d2"-"1d2'" "1d3"-"1d3'"
"1d4"-"1d4'" "1d5"-"1d5'" "1d6"-"1d6'"
"1d7"-"1d7'" "1d8"-"1d8'" "1d9"-"1d9'"
"1u1"-"1u1'" "1u2"-"1u2'" "1u3"-"1u3'" "1u4"-"1u4'"
"1u5"-"1u5'" "1u6"-"1u6'" "1u7"-"1u7'" "1u8"-"1u8'" 
"1u9"-"1u9'" "1u10"-"1u10'" "1u11"-"1u11'" "1u12"-"1u12'"
"2d1"-"2d1'" "2d2"-"2d2'" "2d3"-"2d3'"
"2d4"-"2d4'" "2d5"-"2d5'" "2d6"-"2d6'"
"2d7"-"2d7'" "2d8"-"2d8'" "2d9"-"2d9'"
"2u1"-"2u1'" "2u2"-"2u2'" "2u3"-"2u3'" "2u4"-"2u4'"
"2u5"-"2u5'" "2u6"-"2u6'" "2u7"-"2u7'" "2u8"-"2u8'" 
}\] 
So we have $A_H=\begin{pmatrix}6&1\\4&6 \end{pmatrix}, d=3, \lambda_1=8$, and $\lambda_2=4=\lambda_1^{(d-1)/d}$.
Denote by $T_i^{(m)}, i=1,2$ a tile of type $i$ in $\mathscr{T}^{(m)}$. For every $m\in\N$ we define a patch $V_m\in\mathscr{V}$ in the following process:
\begin{itemize}
\item
Take a tile $T_2^{m+1}$ and remove from it the (unique) $T_1^{(m)}$ that it contains.
\item
From what is left $U_1^{(1)}$, remove all the $T_1^{(m-1)}$ with at least two faces common with $\partial U_1^{(1)}$.
\item[\vdots] 
\item
Eventually, from $U_1^{(m-1)}$ remove all the $T_1^{(1)}$ with at least two faces common with $\partial U_1^{(m-1)}$, to get $U_1^{(m)}$. Define $V_m=U_1^{(m)}$.
\end{itemize}
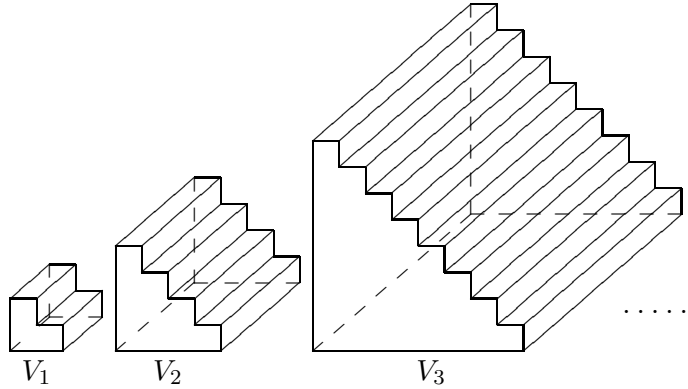
\[\xygraph{
!{<0cm,0cm>;<0.35cm,0cm>:<0cm,0.35cm>::}
!{(0,0)}*{}="A1"
!{(2,0)}*{}="A2"
!{(0,1)}*{}="A3"
!{(1,1)}*{}="A4"
!{(2,1)}*{}="A5"
!{(0,2)}*{}="A6"
!{(1,2)}*{}="A7"
!{(1.5,1.3)}*{}="A'1"
!{(3.5,1.3)}*{}="A'2"
!{(1.5,2.3)}*{}="A'3"
!{(2.5,2.3)}*{}="A'4"
!{(3.5,2.3)}*{}="A'5"
!{(1.5,3.3)}*{}="A'6"
!{(2.5,3.3)}*{}="A'7"
!{(4,0)}*{}="B1"
!{(8,0)}*{}="B2"
!{(7,1)}*{}="B3"
!{(8,1)}*{}="B4"
!{(4,2)}*{}="B5"
!{(6,2)}*{}="B6"
!{(7,2)}*{}="B7"
!{(5,3)}*{}="B8"
!{(6,3)}*{}="B9"
!{(4,4)}*{}="B10"
!{(5,4)}*{}="B11"
!{(7,2.6)}*{}="B'1"
!{(11,2.6)}*{}="B'2"
!{(10,3.6)}*{}="B'3"
!{(11,3.6)}*{}="B'4"
!{(7,4.6)}*{}="B'5"
!{(9,4.6)}*{}="B'6"
!{(10,4.6)}*{}="B'7"
!{(8,5.6)}*{}="B'8"
!{(9,5.6)}*{}="B'9"
!{(7,6.6)}*{}="B'10"
!{(8,6.6)}*{}="B'11"
!{(11.5,0)}*{}="C1"
!{(19.5,0)}*{}="C2"
!{(19.5,1)}*{}="C3"
!{(18.5,1)}*{}="C4"
!{(18.5,2)}*{}="C5"
!{(17.5,2)}*{}="C6"
!{(17.5,3)}*{}="C7"
!{(16.5,3)}*{}="C8"
!{(16.5,4)}*{}="C9"
!{(15.5,4)}*{}="C10"
!{(15.5,5)}*{}="C11"
!{(14.5,5)}*{}="C12"
!{(14.5,6)}*{}="C13"
!{(13.5,6)}*{}="C14"
!{(13.5,7)}*{}="C15"
!{(12.5,7)}*{}="C16"
!{(12.5,8)}*{}="C17"
!{(11.5,8)}*{}="C18"
!{(17.5,5.2)}*{}="C'1"
!{(25.5,5.2)}*{}="C'2"
!{(25.5,6.2)}*{}="C'3"
!{(24.5,6.2)}*{}="C'4"
!{(24.5,7.2)}*{}="C'5"
!{(23.5,7.2)}*{}="C'6"
!{(23.5,8.2)}*{}="C'7"
!{(22.5,8.2)}*{}="C'8"
!{(22.5,9.2)}*{}="C'9"
!{(21.5,9.2)}*{}="C'10"
!{(21.5,10.2)}*{}="C'11"
!{(20.5,10.2)}*{}="C'12"
!{(20.5,11.2)}*{}="C'13"
!{(19.5,11.2)}*{}="C'14"
!{(19.5,12.2)}*{}="C'15"
!{(18.5,12.2)}*{}="C'16"
!{(18.5,13.2)}*{}="C'17"
!{(17.5,13.2)}*{}="C'18"
!{(1,-0.8)}*{V_1}="V1"
!{(6,-0.8)}*{V_2}="V2"
!{(16,-0.8)}*{V_3}="V3"
!{(24.5,1.5)}*{.\:.\:.\:.\:.}="dots"
"A1"-"A2" "A1"-"A6" "A6"-"A7" "A7"-"A4" "A4"-"A5" "A5"-"A2"
"A'1"-@{--}"A'2" "A'1"-@{--}"A'6" "A'6"-"A'7" "A'7"-"A'4" "A'4"-"A'5" "A'5"-"A'2"
"A1"-@{--}"A'1" "A2"-"A'2" "A5"-"A'5" "A4"-"A'4" "A7"-"A'7" "A6"-"A'6"
"B1"-"B2" "B2"-"B4" "B4"-"B3" "B3"-"B7" "B7"-"B6" "B6"-"B9" "B9"-"B8" "B8"-"B11" "B11"-"B10" "B10"-"B1"
"B'1"-@{--}"B'2" "B'2"-"B'4" "B'4"-"B'3" "B'3"-"B'7" "B'7"-"B'6" "B'6"-"B'9" "B'9"-"B'8" "B'8"-"B'11" "B'11"-"B'10" "B'10"-@{--}"B'1"
"B1"-@{--}"B'1" "B2"-"B'2" "B3"-"B'3" "B4"-"B'4" "B6"-"B'6" "B7"-"B'7" "B8"-"B'8" "B9"-"B'9" "B10"-"B'10" "B11"-"B'11"
"C1"-"C2" "C2"-"C3" "C3"-"C4" "C4"-"C5" "C5"-"C6" "C6"-"C7" "C7"-"C8" "C8"-"C9"
"C9"-"C10" "C10"-"C11" "C11"-"C12" "C12"-"C13" "C13"-"C14" "C14"-"C15" "C15"-"C16" "C16"-"C17" "C17"-"C18" "C18"-"C1"
"C'1"-@{--}"C'2" "C'2"-"C'3" "C'3"-"C'4" "C'4"-"C'5" "C'5"-"C'6" "C'6"-"C'7" "C'7"-"C'8" "C'8"-"C'9" "C'9"-"C'10" "C'10"-"C'11" "C'11"-"C'12" "C'12"-"C'13" "C'13"-"C'14" "C'14"-"C'15" "C'15"-"C'16" "C'16"-"C'17" "C'17"-"C'18" "C'18"-@{--}"C'1"
"C1"-@{--}"C'1" "C2"-"C'2" "C3"-"C'3" "C4"-"C'4" "C5"-"C'5" "C6"-"C'6" "C7"-"C'7" "C8"-"C'8" "C9"-"C'9" "C10"-"C'10" "C11"-"C'11" "C12"-"C'12" "C13"-"C'13" "C14"-"C'14" "C15"-"C'15" "C16"-"C'16" "C17"-"C'17" "C18"-"C'18"
}\] 
So obviously 
\begin{equation}\label{eq:Example_boundary_estimate}
\mu_2(\partial V_m)\le\mu_2(\partial T_2^{(m+1)})=6\cdot 4^m.
\end{equation}
We fix an $m$ and estimate $\absolute{\#(Y\cap V_m)-\alpha\cdot\mu(V_m)}$. For that we consider the following partition of $V_m$ to tiles from different levels $\mathscr{T}^{(k)}$'s: 
\[\begin{split}&
\mathscr{U}_m=\{T\in\mathscr{T}^{(m)}:int(T)\subseteq V_m\}
\\&
\mathscr{U}_{m-1}=\{T\in\mathscr{T}^{(m-1)}:int(T)\subseteq V_m\smallsetminus\bigcup\mathscr{U}_m\}
\\& \vdots
\\ 0\le k<m:\quad &
\mathscr{U}_{k}=\{T\in\mathscr{T}^{(k)}:int(T)\subseteq V_m\smallsetminus\bigcup\left(\mathscr{U}_{k+1}\cup\ldots\cup\mathscr{U}_m\right)\}
\end{split}\]
For $i=1,2$ and $k\in\{0,1,\ldots,m\}$ let $t_{i,k}=\#\{T_i^{(k)}\in\mathscr{U}_k\}$. By the construction, 
\begin{equation}\label{eq:Example_t_(i,k)}
t_{1,k}=0\mbox{ for all }k,\mbox{ and }t_{2,k}=\begin{cases}2\cdot 4^{m-k-1},&k\neq 0\\6\cdot 4^{m-1}&k=0\end{cases}.
\end{equation}  
Recall that the discrepancy of $V_m$ depends only on the vector $a_{V_m}=\begin{pmatrix}a_1\\a_2\end{pmatrix}$ (see (\ref{eq:a_V})). We can write it now in terms of the $t_{2,k}$'s. Calculations of $A_H^ke_2$ shows that:
\begin{equation}\label{eq:Example_a_i's}
\begin{split}&
a_1=\sum_{k=0}^mt_{2,k}\cdot A_H^ke_2(1)= 
\sum_{k=0}^m\frac{1}{4}\cdot t_{2,k}(8^k-4^k),
\\&
a_2=\sum_{k=0}^mt_{2,k}\cdot A_H^ke_2(2)= 
\sum_{k=0}^m\frac{1}{2}\cdot t_{2,k}(8^k+4^k)
\end{split}
\end{equation}
Note that $\alpha=3/4$ (see (\ref{eq:alpha})), then
\begin{equation}\label{eq:Example_disc_calculation}
\begin{split}&
\absolute{\#(Y\cap V_m)-\alpha\cdot\mu(V_m)}=
\absolute{a_1+a_2-\frac{3}{4}(2a_1+a_2)}= \absolute{\frac{1}{4}a_2-\frac{1}{2}a_1}
\\&
\stackrel{(\ref{eq:Example_a_i's})}=
\frac{1}{4}\absolute{\sum_{k=0}^mt_{2,k}\cdot 4^k}\stackrel{(\ref{eq:Example_t_(i,k)})}=
\frac{1}{4}\absolute{6\cdot 4^{m-1}+\sum_{k=1}^m2\cdot 4^{m-k-1}\cdot 4^k}=
\left(\frac{m+3}{8}\right)4^m.
\end{split}
\end{equation}
Observe that (\ref{eq:Example_boundary_estimate}) and (\ref{eq:Example_disc_calculation}) together shows that $(ii)$ of Theorem \ref{thm:Laczkovich_Thm} does not hold, which implies that any tiling in $X_H$ correspond to a separated net which is not a BD of $\Z^3$.
\end{proof}

\section{Examples}\label{sec:Examples}
In this last section we give some examples for primitive substitution tilings to show that the different cases that appears in Theorem \ref{thm:Goal} exists. In all of the examples below we give the substitution $H$ and refer the result to any separated net that corresponds to any substitution tiling in $X_H$. Note that in all the examples below the order of the tiles does not matter, but only how many we have of each type. We add the drawings of the substitution rule in order to show that there are substitutions that correspond to the matrices.

\begin{example}\label{example1}
\[\xygraph{
!{<0cm,0cm>;<0.5cm,0cm>:<0cm,0.5cm>::}
!{(0,0)}*{}="a00"
!{(3,0)}*{}="a03"
!{(0,1)}*{}="a10"
!{(1,1)}*{}="a11"
!{(3,1)}*{}="a13"
!{(0,2)}*{}="a20"
!{(1,2)}*{}="a21"
!{(3,2)}*{}="a23"
!{(0,3)}*{}="a30"
!{(3,3)}*{}="a33"
!{(4,0)}*{}="b00"
!{(7,0)}*{}="b03"
!{(10,0)}*{}="b06"
!{(4,1)}*{}="b10"
!{(5,1)}*{}="b11"
!{(7,1)}*{}="b13"
!{(9,1)}*{}="b15"
!{(10,1)}*{}="b16"
!{(4,2)}*{}="b20"
!{(5,2)}*{}="b21"
!{(7,2)}*{}="b23"
!{(9,2)}*{}="b25"
!{(4,3)}*{}="b30"
!{(7,3)}*{}="b33"
!{(9,3)}*{}="b35"
!{(10,3)}*{}="b36"
!{(11,0)}*{}="c00"
!{(12,0)}*{}="c01"
!{(15,0)}*{}="c04"
!{(18,0)}*{}="c07"
!{(19,0)}*{}="c08"
!{(20,0)}*{}="c09"
!{(12,1)}*{}="c11"
!{(15,1)}*{}="c14"
!{(16,1)}*{}="c15"
!{(18,1)}*{}="c17"
!{(19,1)}*{}="c18"
!{(20,1)}*{}="c19"
!{(12,2)}*{}="c21"
!{(14,2)}*{}="c23"
!{(16,2)}*{}="c25"
!{(17,2)}*{}="c26"
!{(18,2)}*{}="c27"
!{(19,2)}*{}="c28"
!{(11,3)}*{}="c30"
!{(12,3)}*{}="c31"
!{(14,3)}*{}="c33"
!{(17,3)}*{}="c36"
!{(19,3)}*{}="c38"
!{(20,3)}*{}="c39"
!{(21,0)}*{}="d00"
!{(25,0)}*{}="d04"
!{(29,0)}*{}="d08"
!{(33,0)}*{}="d012"
!{(21,1)}*{}="d10"
!{(22,1)}*{}="d11"
!{(25,1)}*{}="d14"
!{(26,1)}*{}="d15"
!{(27,1)}*{}="d16"
!{(29,1)}*{}="d18"
!{(30,1)}*{}="d19"
!{(32,1)}*{}="d111"
!{(33,1)}*{}="d112"
!{(21,2)}*{}="d20"
!{(22,2)}*{}="d21"
!{(23,2)}*{}="d22"
!{(26,2)}*{}="d25"
!{(27,2)}*{}="d26"
!{(30,2)}*{}="d29"
!{(31,2)}*{}="d210"
!{(32,2)}*{}="d211"
!{(33,2)}*{}="d212"
!{(21,3)}*{}="d30"
!{(23,3)}*{}="d32"
!{(27,3)}*{}="d36"
!{(31,3)}*{}="d310"
!{(32,3)}*{}="d311"
!{(33,3)}*{}="d312"
"a00"-"a03" "a03"-"a33" "a33"-"a30" "a30"-"a00"
"a10"-"a13" "a20"-"a23" "a11"-"a21"
"b00"-"b06" "b06"-"b36" "b36"-"b30" "b30"-"b00"
"b10"-"b16" "b20"-"b25"
"b03"-"b13" "b11"-"b21" "b23"-"b33" "b15"-"b35"
"c00"-"c09" "c09"-"c39" "c39"-"c30" "c30"-"c00"
"c11"-"c17" "c18"-"c19" "c21"-"c28"
"c01"-"c31" "c04"-"c14" "c07"-"c27" "c08"-"c38"
"c15"-"c25" "c23"-"c33" "c26"-"c36" 
"d00"-"d012" "d012"-"d312" "d312"-"d30" "d30"-"d00"
"d10"-"d112" "d20"-"d212"
"d04"-"d14" "d08"-"d18" "d11"-"d21" "d15"-"d25"
"d16"-"d36" "d19"-"d29" "d111"-"d311" "d22"-"d32" "d210"-"d310"
}\] 
$A_H=\begin{pmatrix} 1&1&1&5 \\ 1&2&5&2 \\ 2&3&4&1 \\ 0&1&1&6 \end{pmatrix}$. The eigenvalues are: $9,4,1,-1$, and we have $4>9^{1/2}$. But the eigenvector that corresponds to $4$ is in $\mathbb{1}^\perp$, then $\lambda_t=1<9^{1/2}$, and therefore any tiling in $X_H$ give rise to a separated net which is a BD of $\Z^2$.
\end{example}

\begin{example}\label{example2}
\[\xygraph{
!{<0cm,0cm>;<0.5cm,0cm>:<0cm,0.5cm>::}
!{(0,0)}*{}="a00"
!{(1,0)}*{}="a01"
!{(2,0)}*{}="a02"
!{(3,0)}*{}="a03"
!{(0,1)}*{}="a10"
!{(1,1)}*{}="a11"
!{(2,1)}*{}="a12"
!{(3,1)}*{}="a13"
!{(0,2)}*{}="a20"
!{(1,2)}*{}="a21"
!{(2,2)}*{}="a22"
!{(3,2)}*{}="a23"
!{(0,3)}*{}="a30"
!{(3,3)}*{}="a33"
!{(4,0)}*{}="b00"
!{(5,0)}*{}="b01"
!{(9,0)}*{}="b05"
!{(10,0)}*{}="b06"
!{(5,1)}*{}="b11"
!{(6,1)}*{}="b12"
!{(8,1)}*{}="b14"
!{(9,1)}*{}="b15"
!{(10,1)}*{}="b16"
!{(4,2)}*{}="b20"
!{(5,2)}*{}="b21"
!{(6,2)}*{}="b22"
!{(7,2)}*{}="b23"
!{(8,2)}*{}="b24"
!{(9,2)}*{}="b25"
!{(10,2)}*{}="b26"
!{(4,3)}*{}="b30"
!{(5,3)}*{}="b31"
!{(7,3)}*{}="b33"
!{(9,3)}*{}="b35"
!{(10,3)}*{}="b36"
!{(11,0)}*{}="c00"
!{(14,0)}*{}="c03"
!{(15,0)}*{}="c04"
!{(16,0)}*{}="c05"
!{(18,0)}*{}="c07"
!{(20,0)}*{}="c09"
!{(11,1)}*{}="c10"
!{(14,1)}*{}="c13"
!{(15,1)}*{}="c14"
!{(16,1)}*{}="c15"
!{(18,1)}*{}="c17"
!{(19,1)}*{}="c18"
!{(20,1)}*{}="c19"
!{(11,2)}*{}="c20"
!{(14,2)}*{}="c23"
!{(15,2)}*{}="c24"
!{(16,2)}*{}="c25"
!{(17,2)}*{}="c26"
!{(19,2)}*{}="c28"
!{(11,3)}*{}="c30"
!{(14,3)}*{}="c33"
!{(17,3)}*{}="c36"
!{(19,3)}*{}="c38"
!{(20,3)}*{}="c39"
!{(21,0)}*{}="d00"
!{(22,0)}*{}="d01"
!{(23,0)}*{}="d02"
!{(25,0)}*{}="d04"
!{(29,0)}*{}="d08"
!{(33,0)}*{}="d012"
!{(23,1)}*{}="d12"
!{(25,1)}*{}="d14"
!{(27,1)}*{}="d16"
!{(29,1)}*{}="d18"
!{(30,1)}*{}="d19"
!{(31,1)}*{}="d110"
!{(33,1)}*{}="d112"
!{(21,2)}*{}="d20"
!{(22,2)}*{}="d21"
!{(23,2)}*{}="d22"
!{(25,2)}*{}="d24"
!{(27,2)}*{}="d26"
!{(29,2)}*{}="d28"
!{(30,2)}*{}="d29"
!{(31,2)}*{}="d210"
!{(32,2)}*{}="d211"
!{(33,2)}*{}="d212"
!{(21,3)}*{}="d30"
!{(25,3)}*{}="d34"
!{(29,3)}*{}="d38"
!{(30,3)}*{}="d39"
!{(32,3)}*{}="d311"
!{(33,3)}*{}="d312"
"a00"-"a03" "a03"-"a33" "a33"-"a30" "a30"-"a00"
"a10"-"a11" "a12"-"a13" "a20"-"a23"
"a01"-"a21" "a02"-"a22"
"b00"-"b06" "b06"-"b36" "b36"-"b30" "b30"-"b00"
"b11"-"b15" "b20"-"b26"
"b01"-"b31" "b21"-"b22" "b23"-"b33" "b14"-"b24" "b05"-"b35"
"c00"-"c09" "c09"-"c39" "c39"-"c30" "c30"-"c00"
"c10"-"c14" "c15"-"c19" "c20"-"c28"
"c03"-"c13" "c23"-"c33" "c04"-"c24" "c05"-"c25"
"c26"-"c36" "c07"-"c17" "c18"-"c38"
"d00"-"d012" "d012"-"d312" "d312"-"d30" "d30"-"d00"
"d12"-"d112" "d20"-"d212"
"d01"-"d21" "d02"-"d22" "d04"-"d14" "d24"-"d34"
"d16"-"d26" "d08"-"d18" "d28"-"d38" "d19"-"d39" 
"d110"-"d210" "d211"-"d311"
}\] 
$A_H=\begin{pmatrix} 4&3&1&3 \\ 1&4&5&5 \\ 1&1&4&1 \\ 0&1&1&5 \end{pmatrix}$. The eigenvalues are: $9,3,3,2$, where $3$ has a non-trivial Jordan block of size $2$, but the eigenvector of $3$ is in $\mathbb{1}^\perp$ and therefore $k_t=1$ (the generalized eigenvector of $3$ is not in $\mathbb{1}^\perp$). So here we cannot apply section $(III)$ of Theorem \ref{thm:Goal} and cannot determine if the corresponding separated nets are BD to $\Z^2$ or not.
\end{example}

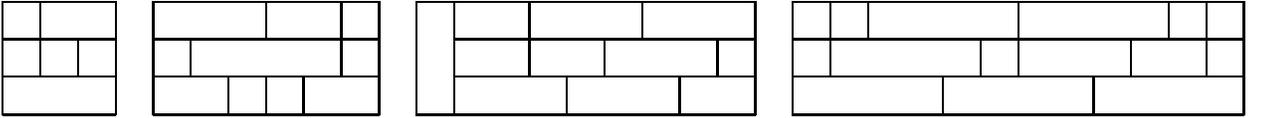
\begin{example}\label{example3}
\[\xygraph{
!{<0cm,0cm>;<0.5cm,0cm>:<0cm,0.5cm>::}
!{(0,0)}*{}="a00"
!{(3,0)}*{}="a03"
!{(0,1)}*{}="a10"
!{(1,1)}*{}="a11"
!{(2,1)}*{}="a12"
!{(3,1)}*{}="a13"
!{(0,2)}*{}="a20"
!{(1,2)}*{}="a21"
!{(2,2)}*{}="a22"
!{(3,2)}*{}="a23"
!{(0,3)}*{}="a30"
!{(1,3)}*{}="a31"
!{(3,3)}*{}="a33"
!{(4,0)}*{}="b00"
!{(6,0)}*{}="b02"
!{(7,0)}*{}="b03"
!{(8,0)}*{}="b04"
!{(10,0)}*{}="b06"
!{(4,1)}*{}="b10"
!{(5,1)}*{}="b11"
!{(6,1)}*{}="b12"
!{(7,1)}*{}="b13"
!{(8,1)}*{}="b14"
!{(9,1)}*{}="b15"
!{(10,1)}*{}="b16"
!{(4,2)}*{}="b20"
!{(5,2)}*{}="b21"
!{(7,2)}*{}="b23"
!{(9,2)}*{}="b25"
!{(10,2)}*{}="b26"
!{(4,3)}*{}="b30"
!{(7,3)}*{}="b33"
!{(9,3)}*{}="b35"
!{(10,3)}*{}="b36"
!{(11,0)}*{}="c00"
!{(12,0)}*{}="c01"
!{(15,0)}*{}="c04"
!{(18,0)}*{}="c07"
!{(19,0)}*{}="c08"
!{(20,0)}*{}="c09"
!{(12,1)}*{}="c11"
!{(14,1)}*{}="c13"
!{(15,1)}*{}="c14"
!{(16,1)}*{}="c15"
!{(18,1)}*{}="c17"
!{(19,1)}*{}="c18"
!{(20,1)}*{}="c19"
!{(12,2)}*{}="c21"
!{(14,2)}*{}="c23"
!{(16,2)}*{}="c25"
!{(17,2)}*{}="c26"
!{(18,2)}*{}="c27"
!{(19,2)}*{}="c28"
!{(20,2)}*{}="c29"
!{(11,3)}*{}="c30"
!{(12,3)}*{}="c31"
!{(14,3)}*{}="c33"
!{(17,3)}*{}="c36"
!{(19,3)}*{}="c38"
!{(20,3)}*{}="c39"
!{(21,0)}*{}="d00"
!{(25,0)}*{}="d04"
!{(29,0)}*{}="d08"
!{(33,0)}*{}="d012"
!{(21,1)}*{}="d10"
!{(22,1)}*{}="d11"
!{(25,1)}*{}="d14"
!{(26,1)}*{}="d15"
!{(27,1)}*{}="d16"
!{(29,1)}*{}="d18"
!{(30,1)}*{}="d19"
!{(32,1)}*{}="d111"
!{(33,1)}*{}="d112"
!{(21,2)}*{}="d20"
!{(22,2)}*{}="d21"
!{(23,2)}*{}="d22"
!{(26,2)}*{}="d25"
!{(27,2)}*{}="d26"
!{(30,2)}*{}="d29"
!{(31,2)}*{}="d210"
!{(32,2)}*{}="d211"
!{(33,2)}*{}="d212"
!{(21,3)}*{}="d30"
!{(22,3)}*{}="d31"
!{(23,3)}*{}="d32"
!{(27,3)}*{}="d36"
!{(31,3)}*{}="d310"
!{(32,3)}*{}="d311"
!{(33,3)}*{}="d312"
"a00"-"a03" "a03"-"a33" "a33"-"a30" "a30"-"a00"
"a10"-"a13" "a20"-"a23" "a11"-"a31" "a12"-"a22"
"b00"-"b06" "b06"-"b36" "b36"-"b30" "b30"-"b00"
"b10"-"b16" "b20"-"b26"
"b02"-"b12" "b03"-"b13" "b04"-"b14" "b11"-"b21" 
"b23"-"b33" "b15"-"b35"
"c00"-"c09" "c09"-"c39" "c39"-"c30" "c30"-"c00"
"c11"-"c19" "c21"-"c29"
"c01"-"c31" "c04"-"c14" "c07"-"c17" "c18"-"c28"
"c15"-"c25" "c13"-"c33" "c26"-"c36" 
"d00"-"d012" "d012"-"d312" "d312"-"d30" "d30"-"d00"
"d10"-"d112" "d20"-"d212"
"d04"-"d14" "d08"-"d18" "d11"-"d31" "d15"-"d25"
"d16"-"d36" "d19"-"d29" "d111"-"d311" "d22"-"d32" "d210"-"d310"
}\] 
$A_H=\begin{pmatrix} 4&5&1&7 \\ 1&3&4&1 \\ 1&1&6&1 \\ 0&1&0&6 \end{pmatrix}$. The eigenvalues are: $9,5,3,2$, and we have $5>9^{1/2}$. But the eigenvector that corresponds to $5$ is in $\mathbb{1}^\perp$, then $\lambda_t=3=9^{1/2}$. Then we have here another example for a substitution that we cannot determine whether the corresponding separated nets are a BD of $\Z^2$ or not.
\end{example}

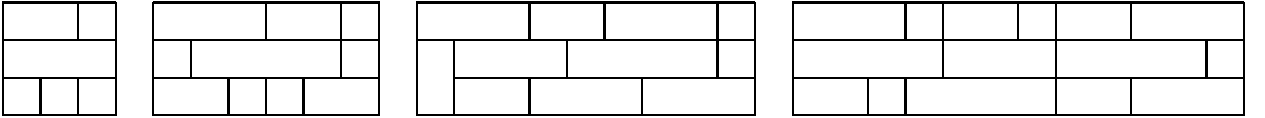
\begin{example}\label{example4}
\[\xygraph{
!{<0cm,0cm>;<0.5cm,0cm>:<0cm,0.5cm>::}
!{(0,0)}*{}="a00"
!{(1,0)}*{}="a01"
!{(2,0)}*{}="a02"
!{(3,0)}*{}="a03"
!{(0,1)}*{}="a10"
!{(1,1)}*{}="a11"
!{(2,1)}*{}="a12"
!{(3,1)}*{}="a13"
!{(0,2)}*{}="a20"
!{(2,2)}*{}="a22"
!{(3,2)}*{}="a23"
!{(0,3)}*{}="a30"
!{(2,3)}*{}="a32"
!{(3,3)}*{}="a33"
!{(4,0)}*{}="b00"
!{(6,0)}*{}="b02"
!{(7,0)}*{}="b03"
!{(8,0)}*{}="b04"
!{(10,0)}*{}="b06"
!{(4,1)}*{}="b10"
!{(5,1)}*{}="b11"
!{(6,1)}*{}="b12"
!{(7,1)}*{}="b13"
!{(8,1)}*{}="b14"
!{(9,1)}*{}="b15"
!{(10,1)}*{}="b16"
!{(4,2)}*{}="b20"
!{(5,2)}*{}="b21"
!{(7,2)}*{}="b23"
!{(9,2)}*{}="b25"
!{(10,2)}*{}="b26"
!{(4,3)}*{}="b30"
!{(7,3)}*{}="b33"
!{(9,3)}*{}="b35"
!{(10,3)}*{}="b36"
!{(11,0)}*{}="c00"
!{(12,0)}*{}="c01"
!{(14,0)}*{}="c03"
!{(17,0)}*{}="c06"
!{(20,0)}*{}="c09"
!{(12,1)}*{}="c11"
!{(14,1)}*{}="c13"
!{(15,1)}*{}="c14"
!{(17,1)}*{}="c16"
!{(19,1)}*{}="c18"
!{(20,1)}*{}="c19"
!{(11,2)}*{}="c20"
!{(12,2)}*{}="c21"
!{(14,2)}*{}="c23"
!{(15,2)}*{}="c24"
!{(16,2)}*{}="c25"
!{(19,2)}*{}="c28"
!{(20,2)}*{}="c29"
!{(11,3)}*{}="c30"
!{(14,3)}*{}="c33"
!{(16,3)}*{}="c35"
!{(19,3)}*{}="c38"
!{(20,3)}*{}="c39"
!{(21,0)}*{}="d00"
!{(23,0)}*{}="d02"
!{(24,0)}*{}="d03"
!{(28,0)}*{}="d07"
!{(30,0)}*{}="d09"
!{(33,0)}*{}="d012"
!{(21,1)}*{}="d10"
!{(23,1)}*{}="d12"
!{(24,1)}*{}="d13"
!{(25,1)}*{}="d14"
!{(30,1)}*{}="d19"
!{(32,1)}*{}="d111"
!{(33,1)}*{}="d112"
!{(21,2)}*{}="d20"
!{(22,2)}*{}="d21"
!{(24,2)}*{}="d23"
!{(25,2)}*{}="d24"
!{(27,2)}*{}="d26"
!{(29,2)}*{}="d28"
!{(30,2)}*{}="d29"
!{(31,2)}*{}="d210"
!{(32,2)}*{}="d211"
!{(33,2)}*{}="d212"
!{(21,3)}*{}="d30"
!{(24,3)}*{}="d33"
!{(25,3)}*{}="d34"
!{(27,3)}*{}="d36"
!{(28,3)}*{}="d37"
!{(30,3)}*{}="d39"
!{(32,3)}*{}="d311"
!{(33,3)}*{}="d312"
"a00"-"a03" "a03"-"a33" "a33"-"a30" "a30"-"a00"
"a10"-"a13" "a20"-"a23" 
"a01"-"a11" "a02"-"a12" "a22"-"a32"
"b00"-"b06" "b06"-"b36" "b36"-"b30" "b30"-"b00"
"b10"-"b16" "b20"-"b26"
"b02"-"b12" "b03"-"b13" "b04"-"b14" "b11"-"b21" "b15"-"b35" "b23"-"b33"
"c00"-"c09" "c09"-"c39" "c39"-"c30" "c30"-"c00"
"c11"-"c19" "c20"-"c29" 
"c01"-"c21" "c03"-"c13" "c06"-"c16" "c14"-"c24" "c18"-"c38"
"c23"-"c33" "c25"-"c35" 
"d00"-"d012" "d012"-"d312" "d312"-"d30" "d30"-"d00"
"d10"-"d112" "d20"-"d212"
"d02"-"d12" "d03"-"d13" "d07"-"d37" "d09"-"d19"
"d14"-"d34" "d111"-"d211" 
"d23"-"d33" "d26"-"d36" "d29"-"d39"
}\] 
$A_H=\begin{pmatrix} 4&5&2&3 \\ 1&3&3&4 \\ 1&1&5&3 \\ 0&1&1&4 \end{pmatrix}$. The eigenvalues are: $9,3,3,1$, where $3$ has a non-trivial Jordan block of size $2$, and both of the vectors are not in $\mathbb{1}^\perp$. So by $(III)$ of Theorem \ref{thm:Goal} the corresponding separated nets are not BD to $\Z^2$.
\end{example}

\end{document}